\documentclass[a4paper]{gtart}

\usepackage[inner=20mm, outer=20mm, textheight=245mm]{geometry}

\usepackage{amsmath}
\usepackage{amssymb}
\usepackage{amsthm,hyperref}

\newtheorem{Theorem}{Theorem}[section]

\newtheorem{Conjecture}[Theorem]{Conjecture}

\newtheorem{Lemma}[Theorem]{Lemma}

\theoremstyle{remark}

\theoremstyle{remark}

\newcommand{\Z}{\mathbb{Z}}
\newcommand{\Q}{\mathbb{Q}}

\DeclareMathOperator{\PGL}{PGL}
\DeclareMathOperator{\PSL}{PSL}
\DeclareMathOperator{\GL}{GL}
\DeclareMathOperator{\SL}{SL}

\usepackage{listings}
\usepackage{color}

\lstdefinestyle{base}{
  emptylines=1,
  moredelim=**[is][\color{blue}]{@}{@},
  moredelim=**[is][\color{green}]{@@}{@@},
}

\usepackage{tikz}
\usepackage{tikz-cd}
\usetikzlibrary{decorations.markings}
\usetikzlibrary{arrows,shapes,positioning}
\usetikzlibrary{matrix}
\usetikzlibrary{cd}
\usetikzlibrary{arrows.meta}

\usepackage{bbm}
\usepackage{mathrsfs}

\begin{document}

\title{Erratum to ``Linear representations of 3-manifold groups over rings''}
\author{Montek Singh Gill\footnotemark[1]}

\begin{abstract}
The fundamental groups of compact $3$-manifolds are known to be residually finite. Feng Luo conjectured that a stronger statement is true, by allowing only finite groups of the form $\PGL_2(R),$ where $R$ is a finite commutative ring with identity. In a prior paper, ``Linear representations of $3$-manifold groups over rings'', we, together with two co-authors, gave an equivalent formulation of Luo's conjecture via faithful representations, investigated variants of the conjecture where the target group is of type $\SL_2, \GL_2, \PSL_2$, or $\PGL_2$ as in the original conjecture, and claimed to provide a $3$-manifold group which violated all four variants of the conjecture. We have discovered a flaw in the work, such that the existing proof is valid only for $\SL_2$ and $\GL_2$. Here, we provide a modified argument for the $\PSL_2$ and $\PGL_2$ cases.
\end{abstract}

\primaryclass{57M27, 57M50}
\keywords{3-manifold, representation of fundamental group, graph manifold}
\makeshorttitle
\footnotetext[1]{Email: montekg@umich.edu}

\section{Introduction}

In \cite{Luo}, Feng Luo made the following conjecture.

\begin{Conjecture}[Luo \cite{Luo}]\label{conj:luo} If $M$ is a compact $3$-manifold and $\gamma \in \pi_1(M) \smallsetminus \{1\}$, then there exists a finite commutative ring $R$ with identity and a homomorphism $\pi_1(M) \rightarrow \PGL_2(R)$ whose kernel does not contain $\gamma$.
\end{Conjecture}

In \cite{FGT}, we considered this conjecture together with the variants where $\PGL_2$ is replaced by any of $\SL_2, \GL_2, \PSL_2$. As in that work, letting $\mathrm{K}$ denote any of $\SL_2, \GL_2, \PSL_2, \PGL_2$, let us refer to the property described in the conjecture as \textit{residual $\mathrm{K}$-finiteness}. The main result of \cite{FGT} is the following.

\begin{Theorem}[Friedl-Gill-Tillmann \cite{FGT}]\label{thm:fgt}
There exists a closed graph manifold $M$ such that $\pi_1(M)$ is not residually $\mathrm{K}$-finite for any of $\mathrm{K} = \SL_2, \GL_2, \PSL_2, \PGL_2$.
\end{Theorem}

We have discovered a flaw in the argument of \cite{FGT}, which breaks the proof in the projective cases $\mathrm{K} = \PSL_2, \PGL_2$. Nevertheless, here, with a modified argument, we demonstrate that the above theorem, with the same 3-manifold group as counterexample, is indeed true.

\section{The Flaw in the Existing Argument}
\label{sec:flaw}

In \cite{FGT}, the main result, Theorem 1.3, is demonstrated in Section 4 (``A Counterexample to Luo's Conjecture''). The 3-manifold group provided is the following:
\[
\Gamma = \langle a,b,c \: | \: ca^2 = b^2, c^{-1}b=bc, ac^{-1}a^{-1} = cac\rangle.
\]
The argument provided relies on the construction of a universal $\SL_2$-representation, together with a SageMath \cite{Sage} computation. The argument does indeed demonstrate the failure of residual $\SL_2$-finiteness, and so of course also of residual $\GL_2$-finiteness. For the projective $\PSL_2, \PGL_2$ cases however, it relies on Proposition 2.9 of \cite{FGT}, and this is where the flaw lies.

\textbf{Erroneous Proposition 2.9 of \cite{FGT}}. \textit{Let $G$ be a finitely generated group. We have the following implications for $G$:}
\begin{center}
 \begin{tikzpicture}
 \draw[-implies, double equal sign distance] (0.2,2.1) node[anchor=east] {\textit{residually $\SL_2$--finite\,\,}} -- (3,2.1)  node[anchor=west] 
 {\textit{\,\,residually $\GL_2$--finite}}
 node[anchor=south,midway] {\textit{}};
 \draw[-implies, double equal sign distance] (3,1.9) -- (0.2,1.9) node[anchor=north,midway] {$Z(G)=1$};
 \draw[implies-, double equal sign distance] (-1.25,1.7) -- (-1.25,0.75) node[anchor=north] {\hspace{-8mm} \textit{residually $\PSL_2$--finite\,\,}} node[anchor=west,midway] {\textit{}};
 \draw[-implies, double equal sign distance] (-1.5,1.7) -- (-1.5,0.75) node[anchor=east,midway] {\textit{$Z(G)$ 2-t.f.}};
 \draw[implies-, double equal sign distance] (4.5,1.7) -- (4.5,0.75) node[anchor = north] 
 {\hspace{13mm} \textit{\hspace{-5mm} \,\,residually $\PGL_2$--finite.}} node[anchor=east,midway] {\textit{}};
 \draw[-implies, double equal sign distance] (4.75,1.7) -- (4.75,0.75) node[anchor=west,midway] {\textit{$\mathrm{Z}(G) = 1$}};
 \draw[-implies,double equal sign distance] (0.4,0.6) -- (3,0.6) node[anchor=south,midway] {\textit{}};
 \draw[-implies,double equal sign distance] (3,0.4) -- (0.4,0.4) node[anchor=north,midway] {\textit{}};
 \end{tikzpicture}
\end{center}
\textit{Here,  2-t.f.\ means 2-torsion-free and $\mathrm{Z}(G)$ denotes the centre of $G$.} \\

The implications which need to be corrected are the unconditional implications which lift from the projective to non-projective cases, that from residual $\PSL_2$-finiteness to residual $\SL_2$-finiteness, and that from residual $\PGL_2$-finiteness to residual $\GL_2$-finiteness, and also the conditional implication from residual $\mathrm{GL}_2$-finiteness to residual $\mathrm{SL}_2$-finiteness. In fact, as argued in \cite{FGT}, in the presence of the other implications, the implication from residual $\mathrm{PGL}_2$-finiteness to residual $\mathrm{GL}_2$-finiteness, and that from residual $\mathrm{GL}_2$-finiteness to $\mathrm{SL}_2$-finiteness, follow from the purported implication from residual $\mathrm{PSL}_2$-finiteness to residual $\mathrm{SL}_2$-finiteness. The mistake is in the proof of this latter implication, from residual $\PSL_2$-finiteness to residual $\SL_2$-finiteness. \\

In \cite{FGT}, in the proof of the implication from residual $\PSL_2$-finiteness to residual $\SL_2$-finiteness, given a representation $G \to \PSL_2(R)$, for some commutative ring $R$, it is claimed that one can construct a factorization $G \to \SL_2(R') \to \PSL_2(R)$, for some other commutative ring $R'$ (the aforementioned implication then follows by Proposition 2.8 of \cite{FGT}). This is done as follows. Let
\[
G = \langle g_1,\dots,g_n \: | \: \{r_i\}_{i\in I} \rangle.
\]
Define the polynomial ring
\[
S = \mathbb{Z}[x_{1a},x_{1b},x_{1c},x_{1d},\dots,x_{na},x_{nb},x_{nc},x_{nd},\{\lambda_i\}_{i\in I}],
\]
and then set
\[
p(g_i) = \left(\begin{array}{cc} x_{ia} & x_{ib} \\ x_{ic} & x_{id} \end{array}\right),\quad  p(g_i^{-1}) = \left(\begin{array}{rr} x_{id} & -x_{ib} \\ -x_{ic} & x_{ia} \end{array}\right).
\]
Define also $p(r_i)$ by setting that $p$ be multiplicative and then define the ideal 
\[
I = \left\langle \{\det\,p(g_i)-1\}_i \cup \{\lambda_i^2-1\}_{i} \cup \{(p(r_i)-\lambda_i )_{k,l}\}_{i,k,l}\right\rangle.
\]
At this point, in \cite{FGT}, $R^\prime$ is set to be the quotient $S/I$ and the desired map $G \to \SL_2(R')$ is defined by
\[
g_i \mapsto \left(\begin{array}{cc} \overline{x_{ia}} & \overline{x_{ib}} \\ \overline{x_{ic}} & \overline{x_{id}} \end{array}\right).
\]
It is claimed that this is a homomorphism. This is not the case though, because, by construction, for each $i$, $p(r_i)$ amounts to the scalar matrix with $\lambda_i$ across the diagonal, \textit{not} the identity matrix. \\

The implication from residual $\PSL_2$-finiteness to residual $\SL_2$-finiteness, and as a result also that from residual $\PGL_2$-finiteness to residual $\GL_2$-finiteness, must be dropped. As for the implication from residual $\mathrm{GL}_2$-finiteness to residual $\mathrm{SL}_2$-finiteness, we can adapt the proof, in \cite{FGT}, of the implication from residual $\mathrm{PGL}_2$-finiteness to residual $\mathrm{PSL}_2$-finiteness to develop an alternate conditional implication. \\

We claim that, if the abelianization of $G$ is torsion-free, then residual $\mathrm{GL}_2$-finiteness implies residual $\mathrm{SL}_2$-finiteness. To see this, fix some element $g_0 \in G$. If $g_0 \notin [G, G]$, then $g_0$ is non-trivial in the abelianization of $G$, and so we may appeal to Proposition 3.5 of \cite{FGT}. Suppose then that $g_0 \in [G,G]$. As $G$ is residually $\mathrm{GL}_2$-finite, there is a finite commutative ring $R$ and a homomorphism
\[
\rho_R \colon G \to \mathrm{GL}_2(R)
\]
such that $\rho_R(g_0)$ is non-trivial. Consider its determinant character
\[
\delta = \det \circ\, \rho_R \colon G \longrightarrow R^\times.
\]
We claim that we can construct a finite ring extension $R \hookrightarrow S$, and a character
\[
\chi \colon G \to S^\times
\]
satisfying
\[
\chi(g)^2 = \delta(g)^{-1}
\]
for all $g \in G$. If we can do this, we can then define a representation
\[
\rho_S \colon G \to \mathrm{SL}_2(S)
\]
by setting
\[
\rho_S(g) = \chi(g)\rho_R(g).
\]
This is a homomorphism because $\chi$ is multiplicative, and maps into $\mathrm{SL}_2(S)$, as opposed to $\mathrm{GL}_2(S)$, because, by construction, we have that
\[
\det\rho_S(g) = \chi(g)^2\det\rho_R(g) = \delta(g)^{-1}\delta(g) = 1.
\]
Moreover, as $S^\times$ is abelian, $\chi$ factors through the abelianization $G_{\mathrm{ab}} = G/[G,G]$, and so, as $g_0 \in [G,G]$, $\chi(g_0) = 1$. It follows that $\rho_S(g_0) = \rho_R(g_0)$, and so $\rho_S(g_0)$ is non-trivial. \\

Thus it remains to construct $S$ and $\chi$. We have already noted that $\chi$ must factor through $G_{\mathrm{ab}}$. We shall construct a map $\bar\chi \colon G_{\mathrm{ab}} \to S^\times$, and the desired $\chi$ is then formed by precomposing with the natural map $G \to G_{\mathrm{ab}}$. By assumption, $G_{\mathrm{ab}}$ is a free abelian group, of rank say $r$. Let $a_1, \dots, a_r$ be free generators of $G_{\mathrm{ab}}$. As with $\chi$, $\delta$ too must factor through the abelianization $G_{\mathrm{ab}}$; let $\bar\delta$ denote the induced map from $G_{\mathrm{ab}}$. For $i = 1, \dots, r$, set $u_i = \bar\delta(a_i)$. Now define $S$ by adjoining square roots to $R$ as follows:
\[
S = R[x_1,\ldots,x_r]\,/\,(x_1^2-u_1^{-1},\ldots,x_r^2-u_r^{-1}).
\]
As in the proof of Proposition 2.9 in \cite{FGT}, the obvious map $R \to S$ is injective. Moreover, since $R$ is finite and $S$ is a finite extension of $R$, $S$ too is finite. Now define $\bar\chi$ by setting, for $i = 1, \dots, r$, that
\[
\bar\chi(a_i) = x_i.
\]
It follows that $\bar\chi(a)^2 = \bar\delta(a)^{-1}$ holds for the generators $a_i$, and so for all $a \in G_{\mathrm{ab}}$ as both sides are homomorphisms. The desired relation $\chi(g)^2 = \delta(g)^{-1}$ for all $g \in G$ follows immediately. This completes the proof. \\

All told, the corrected version of Proposition 2.9 in \cite{FGT} is as below.

\textbf{Corrected Proposition 2.9 of \cite{FGT}}. \textit{Let $G$ be a finitely generated group. We have the following implications for $G$:}
\begin{center}
 \begin{tikzpicture}
 \draw[-implies, double equal sign distance] (0.2,2.1) node[anchor=east] {\textit{residually $\SL_2$--finite\,\,}} -- (3,2.1)  node[anchor=west] 
 {\textit{\,\,residually $\GL_2$--finite}};
 \draw[-implies, double equal sign distance] (3,1.9) -- (0.2,1.9) node[anchor=north,midway] {\textit{$G_{\mathrm{ab}}$ t.f.}};

 \draw[-implies, double equal sign distance] (-1.25,1.7) -- (-1.25,0.75) node[anchor=north] {\hspace{-8mm} \textit{residually $\PSL_2$--finite\,\,}} node[anchor=east,midway] {\textit{$Z(G)$ 2-t.f.}};

 \draw[-implies, double equal sign distance] (4.5,1.7) -- (4.5,0.75) node[anchor = north] 
 {\hspace{13mm} \textit{\hspace{-5mm} \,\,residually $\PGL_2$--finite.}} node[anchor=west,midway] {\textit{$\mathrm{Z}(G) = 1$}}; 
 
 \draw[-implies,double equal sign distance] (0.4,0.6) -- (3,0.6) node[anchor=south,midway] {\textit{}};
 \draw[-implies,double equal sign distance] (3,0.4) -- (0.4,0.4) node[anchor=north,midway] {\textit{}};
 \end{tikzpicture}
\end{center}
\textit{Here, t.f.\ means torsion-free, 2-t.f.\ means 2-torsion-free, $G_{\mathrm{ab}}$ denotes the abelianization of $G$, and $\mathrm{Z}(G)$ denotes the centre of $G$.} \\

We now note the following corrections for the results of \cite{FGT}, all downstream of the correction above for Proposition 2.9:
\begin{itemize}
	\item Corollary 2.10 (1) is directly downstream of Proposition 2.9, and so is false.
	\item Theorem 3.1 for the symmetric groups $S_n$ invokes a SageMath \cite{Sage} computation for the group $S_5$. The computation is in fact for residual $\mathrm{SL}_2$-finiteness, and relies on the incorrect version of Proposition 2.9 for residual $\mathrm{PGL}_2$-finiteness. On the other hand, using the presentation
\[
S_4 = \langle a,b,c \mid a^2 = b^2 = c^2 = 1, (ab)^3 = (bc)^3 = 1, (ac)^2 = 1 \rangle
\]
where $a = (12), b = (23), c = (34)$, one can check that the assignments
\[
a \mapsto \left(\begin{array}{cc} 0 & 1 \\ 1 & 0 \end{array}\right)
\qquad
b \mapsto \left(\begin{array}{cc} 1 & 1 \\ 0 & 1 \end{array}\right)
\qquad
c \mapsto \left(\begin{array}{cc} \varepsilon & 1 \\ 1 & \varepsilon \end{array}\right)
\]
define an embedding $S_4 \hookrightarrow \mathrm{SL}_2(\mathbb{F}_2[\varepsilon]/(\varepsilon^2))$. It follows that Theorem 3.1 holds if we replace $\mathrm{PGL}_2$ with $\mathrm{SL}_2$.
	\item Corollary 3.3 for general linear groups is downstream of Theorem 3.1, and so also holds if we replace $\mathrm{PGL}_2$ with $\mathrm{SL}_2$.
	\item Proposition 3.6 for free groups and surface groups is true, with the proof provided, for $\mathrm{PSL}_2$ and $\mathrm{PGL}_2$, but not necessarily for $\mathrm{SL}_2$ and $\mathrm{GL}_2$.
	\item Proposition 3.7 for geometric manifolds is true, with the argument provided, for $\PSL_2$ and $\PGL_2$, but not a priori for $\SL_2$ and $\GL_2$, except in the case of hyperbolic manifolds, where a separate argument is provided for $\SL_2$ and $\GL_2$.
\end{itemize}

\section{The Corrected Argument}
\label{sec:correction}

The $3$-manifold group in question is the following:
\[
\Gamma = \langle a,b,c \: | \: ca^2 = b^2, c^{-1}b=bc, ac^{-1}a^{-1} = cac\rangle.
\]
In \cite{FGT}, the element of this group for which the residual property is checked is $b^4$. The argument of \cite{FGT} shows that, at this element, $\Gamma$ fails to be residually $\SL_2$-finite, and so of course also residually $\GL_2$-finite. Failure of residual $\PSL_2$-finiteness of course implies failure of residual $\PGL_2$-finiteness, and so it remains for us to check that $\Gamma$ is not residually $\PSL_2$-finite. \\

The simplest method of attack is to compute the image of $b^4$ under the universal $\PSL_2$-representation constructed in \cite{FGT}. We've done this but where the base ring for the polynomial ring in the universal representation is taken to be the rational numbers $\Q$. See the appendix. This demonstrates that, for any characteristic zero ring $R$, any representation $\Gamma \to \PSL_2(R)$ must kill $b^4$. The same computation with $\Q$ replaced by $\Z$ unfortunately proved to be too expensive a computation to succeed. As a result, we first had to simplify the computation, as follows. \\

The universal $\PSL_2$-representation is constructed as follows. First, define the matrices
\[
A = \left(
\begin{array}{cc}
x_{00} & x_{01} \\
x_{10} & x_{11}
\end{array}
\right)
\hspace{1cm}
B = \left(
\begin{array}{cc}
y_{00} & y_{01} \\
y_{10} & y_{11}
\end{array}
\right)
\hspace{1cm}
C = \left(
\begin{array}{cc}
z_{00} & z_{01} \\
z_{10} & z_{11}
\end{array}
\right)
\]
where the $x_{ij}, y_{ij}, z_{ij}$ are indeterminates. Let also $\lambda_1, \lambda_2, \lambda_3$ be indeterminates. Then set $S$ to be the polynomial ring
\[
S = \Z[x_{00}, x_{01}, x_{10}, x_{11}, y_{00}, y_{01}, y_{10}, y_{11}, z_{00}, z_{01}, z_{10}, z_{11}, \lambda_1, \lambda_2, \lambda_3].
\]
Set also $I$ to be the ideal generated by the following:
\begin{itemize}
	\item The polynomials $\det A - 1$, $\det B - 1$, $\det C - 1$.
	\item The polynomials $\lambda_1^2-1$, $\lambda_2^2-1$, $\lambda_3^2-1$.
	\item The polynomial entries of $CA^2 - \lambda_1B^2$, $\mathrm{adj}(C)B - \lambda_2BC$, $A\,\mathrm{adj}(C)\,\mathrm{adj}(A) - \lambda_3CAC$.
\end{itemize}
The universal representation then is given by
\[
\Gamma \to \PSL_2(S/I)
\]
where $a \mapsto [A]$, $b \mapsto [B]$, $c \mapsto [C]$. \\

The following lemma allows us to do away with two of the three scalars $\lambda_i$.

\begin{Lemma}
Let $R$ be a commutative ring and $\rho \colon \Gamma \to \PSL_2(R)$ a representation of $\Gamma$. There exist lifts $A, B, C \in \SL_2(R)$ of $\rho(a)$, $\rho(b)$, $\rho(c)$, respectively, such that
\[
CA^2 = B^2
\hspace{1cm}
C^{-1}B=\mu BC
\hspace{1cm}
AC^{-1}A^{-1} = CAC
\]
for some scalar $\mu \in R$ with $\mu^2 = 1$.
\end{Lemma}

\begin{proof}
As $\rho$ is a representation, there certainly exist lifts $A, B, C$ and scalars $\mu_1, \mu_2, \mu_3$ such that
\[
CA^2 = \mu_1B^2
\hspace{1cm}
C^{-1}B=\mu_2BC
\hspace{1cm}
AC^{-1}A^{-1} = \mu_3CAC
\]
and $\mu_1^2 = \mu_2^2 = \mu_3^2 = 1$. The desired lifts are then $A' = \mu_1\mu_3 A$, $B' = B$, $C' = \mu_1 C$.
\end{proof}

It follows that we can construct a simpler universal representation as follows. Define the matrices $A$, $B$, $C$ as above. Set $S'$ to be the polynomial ring
\[
S' = \Z[x_{00}, x_{01}, x_{10}, x_{11}, y_{00}, y_{01}, y_{10}, y_{11}, z_{00}, z_{01}, z_{10}, z_{11}, \lambda].
\]
Set $I'$ to be the ideal generated by the following:
\begin{itemize}
	\item The polynomials $\det A - 1$, $\det B - 1$, $\det C - 1$.
	\item The polynomial $\lambda^2-1$.
	\item The polynomial entries of $CA^2 - B^2$, $\mathrm{adj}(C)B - \lambda BC$, $A\,\mathrm{adj}(C)\,\mathrm{adj}(A) - CAC$.
\end{itemize}
The universal representation is defined exactly as before. \\

We have now reduced the number of variables in the polynomial ring from $15$ to $13$. We need another reduction. Our goal is to verify that, modulo $I$, $B^4$ is a scalar matrix. Setting $t = \mathrm{tr}\,B$, $d = \det B$, by Cayley-Hamilton, we have that $B^2-tB+d = 0$. From this, it follows that
\[
B^4 = (t^3-2td)B + (d^2-t^2d)I_2
\]
where $I_2$ is the $2 \times 2$ identity matrix. It follows that it suffices to verify that $(t^3-2td)B$ is a scalar matrix. This reduces the verification to ideal membership checks for three degree $4$ polynomials. See now the appendix for this final verification. Note that, even with the above reductions, the computation of the Gr\"{o}bner basis takes non-trivial time (approximately 30 mins on our machine).

\section*{Tool Use}

We used SageMath \cite{Sage} to perform the ideal membership verifications. We used Google's Gemini 4 Argon \cite{Gem}, a large language model, to help make the SageMath code more efficient.

\section*{Appendix}

\subsection*{Ideal Membership Verifications Over $\Q$}

The following is SageMath \cite{Sage} code verifying the claim in Section \ref{sec:correction} that the image of $b^4$ under the universal $\PSL_2$-representation of $\Gamma$, where the base ring is taken to be $\Q$, is trivial.

\begin{lstlisting}[style=base, mathescape=true, escapeinside={<!}{!>}]
S = PolynomialRing(QQ, names=('a00','a01','a10','a11',
                                 'b00','b01','b10','b11',
                                 'c00','c01','c10','c11',
                                 'lmb1', 'lmb2', 'lmb3'),
                     order='degrevlex')
z = S.gens()
A = matrix(S, 2, 2, list(z[0:4]))
B = matrix(S, 2, 2, list(z[4:8]))
C = matrix(S, 2, 2, list(z[8:12]))
lmb1 = z[12]
lmb2 = z[13]
lmb3 = z[14]

def entries(X):
    return [X[0,0], X[0,1], 
             X[1,0], X[1,1]]

def adj(X):
    return matrix(S, 2, 2, [ X[1,1], -X[0,1],
                               -X[1,0],  X[0,0] ])

R1 = C*A**2 - lmb1*B**2
R2 = adj(C)*B - lmb2*B*C
R3 = A*adj(C)*adj(A) - lmb3*C*A*C
relations = ([A.det()-1, B.det()-1, C.det()-1,
               lmb1**2-1, lmb2**2-1, lmb3**2-1]
               + entries(R1) + entries(R2) + entries(R3))
I = S.ideal(relations)

GB = I.groebner_basis()

T = B**4
targets = [T[0,1], T[1,0], T[0,0]-T[1,1]]
answer = tuple(f.reduce(GB) == 0 for f in targets)
print('Scalarity test for b^4 over Q:', answer)

	@Scalarity test for b^4 over Q: (True, True, True)@
\end{lstlisting}

\subsection*{Ideal Membership Verifications Over $\Z$}

The following is SageMath \cite{Sage} code verifying the claim in Section \ref{sec:correction} that the image of $b^4$ under the universal $\PSL_2$-representation of $\Gamma$, where the base ring is taken to be $\Z$, is trivial.

\begin{lstlisting}[style=base, mathescape = true, escapeinside = {<!}{!>}]
S = PolynomialRing(ZZ, names=('a00','a01','a10','a11',
                                 'b00','b01','b10','b11',
                                 'c00','c01','c10','c11',
                                 'lmb'),
                     order='degrevlex')
z = S.gens()
A = matrix(S, 2, 2, list(z[0:4]))
B = matrix(S, 2, 2, list(z[4:8]))
C = matrix(S, 2, 2, list(z[8:12]))
lmb = z[12]

def entries(X):
    return [X[0,0], X[0,1], 
             X[1,0], X[1,1]]

def adj(X):
    return matrix(S, 2, 2, [ X[1,1], -X[0,1],
                               -X[1,0],  X[0,0] ])

R1 = C*A**2 - B**2
R2 = adj(C)*B - lmb*B*C
R3 = A*adj(C)*adj(A) - C*A*C
relations = ([A.det()-1, B.det()-1, C.det()-1,
               lmb**2-1]
               + entries(R1) + entries(R2) + entries(R3))
I = S.ideal(relations)

GB = I.groebner_basis()

t = B.trace()
d = B.det()
h = t*(t*t-2*d)
targets = [h*B[0,1], h*B[1,0], h*(B[0,0]-B[1,1])]
answer = tuple(f.reduce(GB) == 0 for f in targets)
print('Scalarity test for b^4 over Z:', answer)

	@Scalarity test for b^4 over Z: (True, True, True)@
\end{lstlisting}

\end{document}